\documentclass[12pt]{amsart}

\usepackage[english]{babel}
\usepackage{amsfonts}
\usepackage{amscd}
\usepackage{amsmath}

\newtheorem{Thm}{Theorem}
\newtheorem{rmk}{Remark}
\newenvironment{Rmk}{\begin{rmk}\em}{\end{rmk}}
\newtheorem{exm}{Example}

\newtheorem{prf}{Proof}

\newenvironment{Prf}{\begin{prf}\em}{\qed\end{prf}}
\newtheorem{prff}{}

\newtheorem{stat}{}

\newcommand{\NOT}[1]{}

\newcommand{\pa}{\par\medskip}

\newcommand{\LP}{\left(}  \newcommand{\RP}{\right)}

\newcommand{\BE}{\begin{equation}}  \newcommand{\EE}{\end{equation}}
\newcommand{\BR}{\begin{eqnarray*}}  \newcommand{\ER}{\end{eqnarray*}}
\newcommand{\BER}{$$\begin{array}}  \newcommand{\EER}{\end{array}$$}

\newcommand{\RRe}{{\text{Re}\,}}
\newcommand{\IIm}{{\text{Im}\,}}

\newcommand{\I}{\infty}

\newcommand{\al}{\alpha}
\newcommand{\be}{\beta}
 \newcommand{\ga}{\gamma}

 \newcommand{\la}{\lambda}

\newcommand{\bbb}[1]{\mathbb{#1}}

\newcommand{\bC}{\bbb{C}}

\newcommand{\cQ}{\mathcal{Q}}

\newcommand{\tr}{{\text{tr}\,}}

\newcommand{\U}{\mathbf{I}}

\newcommand{\FQ}{\LP\frac{\cQ(A_0)}{|\cQ(A_0)|}\RP^{-1/2}}

\title{A Bound Below for the Convex Hull of the Spectrum of a Matrix}

\author{Eliahu Levy}
\address{Department of Mathematics,
Technion -- Israel Institute of Technology,
Haifa 32000, Israel}
\email{eliahu@techunix.technion.ac.il}

\date{}

\begin{document}

\begin{abstract}
In this note the following is shown. Consider the quadratic form on (complex) matrices $\cQ(A):=\tr(A^2)$. Let $A$ be such a matrix. Then an ellipse can be found, with the vector from center to focus determined by the value of $\cQ$ at the traceless part of
$A$, which must be contained in the convex hull of the spectrum of $A$.
\end{abstract}

\maketitle

Denote by $\bC$ the set of complex numbers, and by $\bC^{n\times n}$ the set of $n\times n$ matrices over $\bC$. Denote by $\U$ the unit matrix. \pa

Consider the quadratic form $\cQ(A):=\tr(A^2)$ = the sum of squares of the eigenvalues of $A$ (which we always count according to their multiplicities as roots of the characteristic polynomial), $A\in\bC^{n\times n}$. ($\cQ$ is a quadratic, not a Hermitian form!) \pa

The symmetric bilinear form (``inner product'') corresponding to $\cQ$ is $A,B\mapsto\tr(AB)=\tr(BA)$. \pa

Note that $\cQ$ is {\em invariant under similarity of matrices}. Thus, $\cQ$ is preserved by any similarity transformation of $\bC^{n\times n}$\,\,$A\mapsto T^{-1}AT$,\,\,$T$ fixed invertible. \pa

Clearly, $\cQ$ is nondegenerate, i.e.\ only for the matrix $A=0$ we have $\tr(AX)=0$ for all $X\in\bC^{n\times n}$.
\par\bigskip

We wish to say something about the convex hull of the spectrum (hence about the spectral radius) of a matrix $A\in\bC^{n\times n}$, decomposed as
\BE\label{eq:dec}
A=\ga\U+A_0,\qquad\ga=\dfrac1n\tr A\in\bC,\quad\tr A_0=0.
\EE
Note that the ``inner product'' $\tr(A_0\cdot\ga\U)=0$. Therefore
\BE\label{eq:dec1}
\cQ(A)=n\ga^2+\cQ(A_0)=\frac1n(\tr A)^2+\cQ(A_0).
\EE

\begin{Thm}\label{Thm:ellipse}
Let $A_0\in\bC^{n\times n}$ be traceless. Then the convex hull of the spectrum of $A_0$ contains an ellipse with foci $\pm\frac1{\sqrt2(n-1)}\sqrt{\cQ(A_0)}$ and with sum of squares of the semiaxes equal to $\LP\frac1{\sqrt2(n-1)}\RP^2$ times the sum of squares of the absolute values of the eigenvalues of $A_0$ (counted by multiplicity in the characteristic polynomial). Indeed, if $\cQ(A_0)\ne0$ then the semiaxes of this ellipse are: the semimajor axis
$$\frac1{\sqrt2(n-1)}
\sqrt{\LP\RRe\LP\la_1\FQ\RP\RP^2\NOT{+\LP\RRe\LP\la_2\FQ\RP\RP^2}+\ldots+
\LP\RRe\LP\la_n\FQ\RP\RP^2},$$
and the semiminor axis
$$\frac1{\sqrt2(n-1)}
\sqrt{\LP\IIm\LP\la_1\FQ\RP\RP^2\NOT{+\LP\IIm\LP\la_2\FQ\RP\RP^2}+\ldots
+\LP\IIm\LP\la_n\FQ\RP\RP^2},$$
where $\la_1,\ldots,\la_n$ are the eigenvalues of $A_0$, counted by multiplicity.

Consequently, for a general $A$, the convex hull of the spectrum will include the shifted
ellipse with foci $\ga\pm\frac1{\sqrt2(n-1)}\sqrt{\cQ(A_0)}$, with $\ga$ and $A_0$ as in
(\ref{eq:dec}), and, of course, will contain the analogous ellipse for any invariant subspace with respect to $A$.
\end{Thm}
\begin{Prf}
Let the eigenvalues of $A_0$ be $\la_1,\ldots,\la_n$.
Then (note $A_0$ is traceless):
\BR
\la_1+\la_2+\ldots+\la_n&=&0,\\
\la_1^2+\la_2^2+\ldots+\la_n^2&=&\cQ(A_0).
\ER
In case $\cQ(A_0)\ne0$ define
\BE\label{eq:mula}
\mu_i:=\la_i\FQ,\quad i=1,\ldots,n
\EE
(fix one of the square roots). If $\cQ(A_0)=0$ let $\mu_i:=\la_i$,\,\,$i=1,\ldots,n$. Then
\BR
\mu_1+\mu_2+\ldots+\mu_n&=&0,\\
\mu_1^2+\mu_2^2+\ldots+\mu_n^2&=&|\cQ(A_0)|.
\ER
Taking real and imaginary parts, this gives, denoting
\BR
R&:=&\sqrt{(\RRe\mu_1)^2+(\RRe\mu_2)^2+\ldots+(\RRe\mu_n)^2},\\
I&:=&\sqrt{(\IIm\mu_1)^2+(\IIm\mu_2)^2+\ldots+(\IIm\mu_n)^2},
\ER
that
\BR
\RRe\mu_1+\RRe\mu_2+\ldots+\RRe\mu_n&=&0,\\
\IIm\mu_1+\IIm\mu_2+\ldots+\IIm\mu_n&=&0,\\
(\RRe\mu_1)(\IIm\mu_1)+(\RRe\mu_2)(\IIm\mu_2)+\ldots+(\RRe\mu_n)(\IIm\mu_n)&=&0,\\
R^2-I^2&=&|\cQ(A_0)|.
\ER
Hence, for any fixed real $\al$ and $\be$,
\begin{eqnarray}
&&(\al\RRe\mu_1+\be\IIm\mu_1)+(\al\RRe\mu_2+\be\IIm\mu_2)+\ldots+
(\al\RRe\mu_n+\be\IIm\mu_n)=0,\label{eq:eq1}\\
&&(\al\RRe\mu_1+\be\IIm\mu_1)^2+(\al\RRe\mu_2+\be\IIm\mu_2)^2+\ldots+
(\al\RRe\mu_n+\be\IIm\mu_n)^2=\al^2R^2+\be^2I^2.\label{eq:eq2}
\end{eqnarray}
Fixing $\al$ and $\be$, suppose there are $n_p$ nonnegative $\al\RRe\mu_i+\be\IIm\mu_i$'s, forming an $n_p$-tuple $\vec{v}_p$, and $n_n$ negative $\al\RRe\mu_i+\be\IIm\mu_i$'s, forming an $n_n$-tuple $\vec{v}_n$, $n_p+n_n=n$. By (\ref{eq:eq1}) the sum of entries of the positive vector $\vec{v}_p$ is minus the sum of entries of the negative $\vec{v}_n$, i.e.\ these vectors have the same $\ell_1$-norm $L$. By (\ref{eq:eq2}) the sum of the squares of the $\ell_2$-norms of $\vec{v}_p$ and $\vec{v}_n$ is $\al^2R^2+\be^2I^2$, and since one always has $|\;|_2^2\le|\;|_1\cdot|\;|_\I$, we have
$$\al^2R^2+\be^2I^2=|\vec{v}_p|_2^2+|\vec{v}_n|_2^2\le
L\cdot(|\vec{v}_p|_\I+|\vec{v}_n|_\I).$$
Also $|\vec{v}_p|_\I\le L\le(n-1)|\vec{v}_p|_\I$,\,\,
$|\vec{v}_n|_\I\le L\le(n-1)|\vec{v}_n|_\I$.
Thus,
$$\al^2R+\be^2I\le2L^2,$$
$$L\ge\sqrt{\frac12(\al^2R^2+\be^2I^2)},$$
$$|\vec{v}_p|_\I,\;|\vec{v}_n|_\I\ge\frac{L}{n-1}\ge
\frac1{\sqrt2(n-1)}\sqrt{\al^2R^2+\be^2I^2}.$$
In particular, this means that
\BE\label{eq:con}
\max_i(\al\RRe\mu_i+\be\IIm\mu_i)\ge\frac1{\sqrt2(n-1)}\sqrt{\al^2R^2+\be^2I^2}.
\EE
Now, $\sqrt{\al^2R^2+\be^2I^2}$ is the maximum of $\al R\cdot\RRe\zeta'+\be  I\cdot\IIm\zeta'$ for $\zeta'$ in the unit disk. Therefore it is the maximum of $\al\RRe\zeta+\be\IIm\zeta$ for $\zeta$ in the ellipse with axes $R$ at the real axis and $I$ at the imaginary axis, which, since $R^2-I^2=|\cQ(A_0|$, has $\pm\sqrt{|\cQ(A_0)|}$ as foci. By convex separation in the plane, (\ref{eq:con}) holding for any real $\al$ and $\be$ implies that the convex hull of $\{\mu_1,\mu_2,\ldots,\mu_n\}$ contains that ellipse multiplied by $\frac1{\sqrt2(n-1)}$. Using (\ref{eq:mula}) to transfer that to working with the $\la_i$'s, one has our assertion.
\end{Prf}

\begin{Rmk}
Compare Theorem \ref{Thm:ellipse} with the example of a traceless $A$ with $n-1$ eigenvalues $-1$ and one eigenvalue $n-1$. Here $\sqrt{\cQ(A)}=\sqrt{n(n-1)}$ and the convex hull of the spectrum is $[-1,n-1]$, so only a major axis $\le1$ for the ellipse may do. Thus one cannot avoid the $n$ in the denominator in Theorem \ref{Thm:ellipse}.

\end{Rmk}


\end{document}